\documentclass[10pt,english]{amsart}

\usepackage{helvet}
\usepackage[T1]{fontenc}
\usepackage[latin1]{inputenc}
\usepackage{a4wide}
\usepackage{amssymb}
\usepackage{babel}

\usepackage{enumerate}
\usepackage{amsfonts}
\usepackage{amsmath}
\usepackage{amsthm}
\usepackage{amssymb}
\usepackage{latexsym}
\usepackage{multicol}
\usepackage{verbatim}
\usepackage{amscd}
\usepackage{hyperref}
\usepackage[usenames]{color}
\usepackage{pb-diagram}

\newtheorem{theorem}{Theorem}
\newtheorem{definition}[theorem]{Definition}
\newtheorem{corollary}[theorem]{Corollary}
\newtheorem{lemma}[theorem]{Lemma}
\newtheorem{remark}[theorem]{Remark}
\newtheorem{example}[theorem]{Example}

\newcommand{\trdeg}{\operatorname{trdeg}}
\newcommand{\aut}{\operatorname{Aut}}
\newcommand{\der}{\operatorname{Der}}
\newcommand{\lnd}{\operatorname{LND}}
\newcommand{\klnd}{\operatorname{KLND}}
\newcommand{\Frac}{\operatorname{Frac}}

\begin{document}

\title{Locally Nilpotent Derivations and Automorphism groups of certain Danielewski surfaces}

\author{Angelo Calil Bianchi and Marcelo Oliveira Veloso }
\address[A.~Bianchi]{UNIFESP - Institute of Science and Technology, S\~ao Jos\'e dos Campos - SP}
\email{acbianchi@unifesp.br}
\address[M.~Veloso]{UFSJ - Department of Physics and Mathematics, Alto Paraopeba - MG}
\email{veloso@ufsj.edu.br}
\date{\today}

\begin{abstract}
We describe the set of all locally nilpotent derivations of the quotient ring \linebreak $\mathbb{K}[X,Y,Z]/(f(X)Y - \varphi(X,Z))$ constructed from the defining equation $f(X)Y = \varphi(X,Z)$ of a \textit{generalized Danielewski surface} in $
\mathbb K^3$ for a specific choice of polynomials $f$ and $\varphi$, with  $\mathbb K$ an algebraically closed field of characteristic zero. As a consequence of this description we calculate the $ML$-invariant and the Derksen invariant of this ring. We also determine a set of generators for the group of $\mathbb K$-automorphisms of $\mathbb K[X,Y,Z]/(f(X)Y - \varphi(Z))$ also for a specific choice of polynomials $f$ and $\varphi$.
\end{abstract}

\maketitle



\section*{Introduction}

The term \textit{Danielewski surface} usually refers to surfaces given by an equation of the form $X^nZ= P(Y)$, with $n\in\mathbb N$ and certain polynomials $P(Y)\in \mathbb C[Y]$, because such surfaces were studied by Danielewski in connection with the famous Cancellation Problem (see \cite{lml3}). Its generalizations continue to be a source of interest for current research. 

For certain authors, a Danielewski surface is an affine surface which is algebraically isomorphic to a surface defined by an equation of the form $X^nZ-P(Y)=0$, for (specific choice of) $P(Y) \in \mathbb C[Y]$, or even a surface defined by the equation $X^nZ-Q(X,Y)=0$, where $Q(X,Y)$is a polynomial satisfying certain properties. Over the past 30 years, many works on surfaces given by this type of equation were published under algebraic and algebraic-geometric approach (see \cite{dan,dda2}). Following we trace some papers in the same direction which we are interested: 
\begin{itemize}
	\item L. Makar-Limanov in \cite{lml1,lml3} computed automorphism groups of surfaces in $\mathbb C^3$ defined by equation with the form $X^nZ-P(Y)=0$, where $n\ge 1$ and $P(Y)$ is a nonzero polynomial. The $ML$-invariant is used in \cite{lml3} to find the group of $\mathbb K$-automorphisms of the ring  $\mathbb K[X,Y,Z]/(X^nZ-P(Y))$, where $n>1$ and $\deg P \ge 2$.

  \item D. Daigle  in \cite{dda1} studied the locally nilpotent derivations of the ring $R=\mathbb K[X,Y,Z]/(XY - \varphi(Z))$ and showed that certain subgroups of $\mathbb K$-automorphisms of $R$ act transitively on the kernels of the nontrivial locally nilpotent derivations on $R$.

  \item A. Crachiola in \cite{cra1} obtained similar results for slightly different surfaces defined by the equations $X^nZ - Y^2 - \sigma(X) Y=0$, where $n\ge 2$ and $\sigma(0)\ne 0$, defined over arbitrary base field.

\item A. Dubouloz and P-M Poloni \cite{dubo} considered more general surfaces defined by equations $X^nZ-Q(X,Y)=0$, where $n\ge 2$ and $Q(X,Y)$ is a polynomial with coefficients in an arbitrary base field such that $Q(0, Y)$ splits with $r\ge 2$ simple roots. This class contains most of the surfaces considered by L. Makar-Limanov, D. Daigle, and A. Crachiola. 
\end{itemize}
  
In this paper we obtain some similar results for a class of Danielewski surfaces given by the equation $f(X)Y - \varphi(X,Z)=0$, that means we study the ring $$\mathcal B=\mathbb K[X,Y,Z]/(f(X)Y - \varphi(X,Z)),$$ where $X,Y$ and $Z$ are indeterminates over $\mathbb K$, $\varphi(X,Z)= Z^m + b_{m-1}(X)Z^{m-1}+ \dots + b_1(X)Z + b_0(X)$, $\mathbb K$ is an algebraically closed field of characteristic zero, $\deg(f)>1$, and $m>1$. We also may write $\mathcal B=\mathbb K[x,y,z]$ where $x,y$ and $z$ are the images of $X,Y$ and $Z$ under the canonical epimorphism $\mathbb K[X,Y,Z] \to \mathcal B$. Note that we have $f(x)y-\varphi(x,z)=0$. 

We show that the $ML$-invariant of $\mathcal B$ is $\mathbb K[x]$ (cf. Theorem \ref{akkx}) and, by using this result, we describe all locally
nilpotent derivations of $\mathcal B$ (cf. Corollary \ref{lndb}) and its Derksen invariant (cf. Theorem \ref{HD}). Further we determine a set of generators for the group of $\mathbb K$-automorphisms of $\mathcal B$ (cf. Theorem \ref{autb}), when $\varphi(X,Z)\in \mathbb K[Z]$ and $f(X)$ has at least one nonzero root which is the case not covered by any of the previously mentioned papers. For all these results we are specially motivated by \cite{dda2, lml3}.

The material is organized as follows: Section 1 gathers the basic definitions, notations, and results used in this paper. In Section 2 we discuss several properties of the ring $\mathcal B$. The set of locally nilpotent $\mathbb K[x]$-derivations of $\mathcal B$ is described in Section 3. The $ML$-invariant and the Derksen invariant of $\mathcal B$ are calculated in Section 4. Finally, Section 5 describes a set of generators for the group of $\mathbb K$-automorphisms of $\mathcal B$ with the assumption that $f(X)$ has at least one nonzero root and $\varphi(X,Z)\in \mathbb K[Z]$.

In a forthcoming paper we will deal with a generalization including more than three variables.

\section{Preliminaries} \label{pre}

\subsection{Some generalities.} From now on, word  ``\textit{ring}'' means a commutative ring with unity and characteristic zero.
We denote the group of units of a ring $A$ by $A^*$. 
In the case that a ring $A$ is a domain, we denote by $\Frac(A)$ its fraction field; when $B$ is a subring of a ring $A$ and $A$ is an integral domain, we denote by $\trdeg_B(A)$ the transcendence degree of $\Frac(A)$ over $\Frac(B)$.
We reserve capital letters $X,Y,Z$ to denote indeterminates over a field $A$ and we denote the polynomial ring in $X,Y,Z$ over $A$ by $A[X,Y,Z]$. The polynomial ring in $n$-indeterminates over $A$ is also denoted by $A^{[n]}$. 
We reserve the symbol $\mathbb K$ to indicate a field.

\subsection{Locally nilpotent derivations.}  Let $A$ be a ring. A derivation $D: A \rightarrow A$ is called \textit{locally nilpotent} if for each $a\in A$ there exist $n\in \mathbb N$ (depending on $a$) such that $D^n(a)=0$.  

In the case of polynomial rings, the basic example of locally nilpotent derivations are the partial derivatives. Given a polynomial $G\in \mathbb K[X_1,\dots,X_n]$, we may write $G_{X_i}$ to denote the partial derivative $\frac{\partial (G)}{\partial {X_i}}$. 

Let us also use the following notations:

\hspace{1cm} $\ker(D)=\{a \in A \mid D(a)=0\}$

\hspace{1cm} $\der(A)=\{D \mid D: A \rightarrow A \:\: \mbox{ is a derivation of } A\}$

\hspace{1cm} $\lnd(A)= \{D \in \der(A) \mid D\neq 0 \mbox{ and } D\mbox{ is locally nilpotent} \}$

\hspace{1cm} $\lnd_R(A)= \{D \in \lnd(A) \mid D(R)=0 \} $

\hspace{1cm} $\klnd(A)= \{ \ker D \mid D \in \lnd(A)\}.$

Following we collect some well known results which will be used throughout the text. For the next three lemmas, details can be found in \cite{dda3,gf1}.

\begin{lemma}\label{tjik}
Let $A$ be a domain and  $D \in \der(A)$.
\begin{enumerate}
\item  If $D \in \lnd(A)$ then $\trdeg_{\ker D}(A)=1$ and $\ker(\Frac(A))=\Frac(\ker A)$.

\item  If $A$ is a subring of $\mathbb K(X_1, \ldots,X_n)$ such that $\trdeg_{\mathbb K}(\ker D)=n-1$ and $f_1,\dots,f_{n-1} \in \ker D$ is a transcendence  basis over $\mathbb K$, then there exists $h \in \mathbb K(X_1, \ldots,X_n)$ such that $D(a)=h J(f_1, \dots, f_{n-1},a)$ for all $a \in A$, where $J(f_1,\dots, f_{n-1},a)$ is the Jacobian relative to $(x_1, \dots,x_n)$.

\item If $D_1,D_2 \in \lnd(A)$, $\ker D_1 = \ker D_2 := C$ and there exists $s \in A$ such that $D_1(s) \in C\setminus\{0\}$, then $D_2(s) \in C\setminus\{0\}$ and $D_2(s)D_1=D_1(s)D_2$.  
\end{enumerate} 
\end{lemma}

\begin{lemma}\label{prop}
Let $A$ be a domain and  $D \in \lnd(A)$.
\begin{enumerate}
\item $\ker D$ is factorially closed in $A$ (that is: if $a,b \in A\setminus \{ 0\}$ and $ab \in \ker D$, then $a, b\in \ker D$).

\item $\ker D$ is algebraically closed in $A$.

\item If $A$ is an unique factorization domain  (UFD), then so is $\ker D$.
\end{enumerate}
\end{lemma}

Recall that a derivation $D\in \der(A)$ is \textit{irreducible} if the unique principal ideal of $A$ which contains $D(A)$ is $A$.

\begin{lemma} \label{clnd} Let $R$ be a domain of characteristic zero satisfying the ascending chain condition for principal ideals, and  $A \in \klnd(R)$. Consider the set $S=\{D \in \lnd_A(R) \mid D \mbox{ is irreducible }\}$. Then $S\neq \emptyset$ and $\lnd_A(R) =\{aD  \mid a \in A \mbox{ and } D \in S\}$. 
\end{lemma}

\subsection{Gradation and filtration.} Let $A$ be a $\mathbb K$-algebra. 

We say that $A$ is $\mathbb{Z}$-graded if there exists a family $\{B_i\}_{i \in \mathbb Z}$ of linear subspaces of $A$ such that $B_iB_j \subseteq B_{i+j}$ for all $i,j \in \mathbb{Z}$ and $A=\displaystyle\oplus_{i \in \mathbb{Z}} B_i$. Similarly, we say that $A$ has a $\mathbb{Z}$-filtration if there exists a family $\{A_i\}_{i \in \mathbb{Z}}$ of linear subspaces of $A$ such that 

\begin{enumerate}
\item[(i)] $A_iA_j\subseteq
A_{i+j}$ for all $i,j \in \mathbb{Z}$,
\item[(ii)] $A_i\subseteq A_j$ for all $i,j \in \mathbb{Z}$ with $i\leq j$, and
\item[(iii)] $A=\cup_{i \in \mathbb{Z}} A_i$.
\end{enumerate}
\noindent If we also have
\begin{enumerate}
\item[(iv)]  $\cap_{i\in \mathbb{Z}}A_i=\{0\}$,
\item[(v)]  $a\in A_n\setminus A_{n-1}$ and $b\in A_m\setminus A_{m-1}$ $\implies$ $ab\in A_{m+n}\setminus A_{m+n-1}$,
\end{enumerate}
we say that $A$ has a proper $\mathbb{Z}$-filtration.

\

Let  $\{A_i\}_{i \in \mathbb{Z}}$ be a $\mathbb{Z}$-filtration of $A$ satisfying (iv) and $D \in \der(A)$ such that $D(A_i)\subseteq A_{i+k}$ for a fixed  $k\in \mathbb N$ and all $i \in {\mathbb Z}$. Let $Gr(A)=\oplus_{i \in \mathbb{Z}} \frac{A_i}{A_{i-1}}$ be the corresponding graded algebra. For $\overline{h} \in  \frac{A_i}{A_{i-1}}$, write $\overline{h}= h + A_{i-1}$ for some $h \in A_i$ and define a derivation $\mathbb{D}$ on $Gr(A)$ by linearly extending the homomorphism defined by $\mathbb{D}(\overline{h}):=D(h)+A_{i+k-1} \in \frac{A_{i+k}}{A_{i+k-1}}$. If $k$ is chosen minimal (but not $-\infty$), then $\mathbb{D} \ne 0$. Defining the map $gr: A \to Gr(A)$ by $gr(a)= a + A_{i-1}$ if $a \in A_i \setminus A_{i-1}$, and $gr(0)=0$, then either $\mathbb{D}(gr(a))=gr(D(a))$ or $\mathbb{D}(gr(a))= 0$. 
Further, if $D\in LND(A)$, then $\mathbb{D}$ is locally nilpotent derivation on $Gr(A)$ (see \cite[Lemma 4]{lml2}).

\section{The main object \texorpdfstring{$\mathcal B$}{}}

Let $\mathcal B=\mathbb K[X,Y,Z]/(f(X)Y - \varphi(X,Z))$, where $\mathbb K$ is an algebraically closed field of characteristic zero, $f(X)=X^r+a_{r-1}X^{r-1}+\dots+a_0\in \mathbb K[X]$ with $\deg_X f = r>1$, and $\varphi(X,Z)=Z^d +b_{d-1}(X)Z^{d-1}+\dots + b_1(X)Z+b_0(X)\in \mathbb K[X][Z]$ with  $\deg_Z \varphi = d>1$.  Notice that we also can write $\mathcal B=\mathbb K[x,y,z]$, where $x,y$ and $z$ are the images of $X,Y$ and $Z$ under the canonical epimorphism $\mathbb K[X,Y,Z] \to \mathcal B$. We now state some basic properties of the ring $\mathcal B$.

\begin{lemma}  \label{gmn}  \
\begin{enumerate}

\item The ring $\mathcal B$ is a domain and $\trdeg_\mathbb K(\mathcal B)=2$. 

\item For each $b \in \mathcal B\setminus \{0\}$, there exists a unique $g \in \mathbb K[X,Y,Z]$ such that $\deg_{Z} (g) < d$ and $b=g(x,y,z)$.

\item $\mathbb K(x)\cap \mathcal B=\mathbb K[x]$.

\item Let $S=\mathbb K[x]\setminus \{0\}$. Then, the localization o $\mathcal B$ by $S$ satisfies $S^{-1}\mathcal B=\mathbb K(x)^{[1]}$.

\item $\mathbb K[x]$ is factorially closed in $\mathcal B$ and $\mathcal B^{*}=\mathbb K^{*}$.
\end{enumerate}

\end{lemma}

\proof To prove (1), note that $F(X,Y,Z)=f(X)Y-\varphi(X,Z)$ is  irreducible in $\mathbb K[X,Y,Z]$. Item (2) follows directly from the division algorithm by considering $F(X,Y,Z)$ as a monic polynomial in $\mathbb K[X,Y][Z]$.
For item (3), let $b\in \mathbb K(x)\cap \mathcal B\setminus \{0\}$ then $b=\sum_{i<d}a_iz^i$, where $a_i \in \mathbb K[x,y]$, and  there exists $a \in \mathbb K[x]\setminus \{0\}$ such that  $ab \in \mathbb K[x]$. Thus $ab= \sum_{i<d}(aa_i)z^i$, which implies that $a_i=0$ for all $i>0$; so $b=a_0 \in \mathbb K[x,y]\cap \mathbb K(x)=\mathbb K[x]$ and the statement is true. To show (4), since $y=f(x)^{-1}\varphi(x,z) \in \mathbb K(x)[z]$, we have $\mathbb K[x,z]\subseteq \mathcal B \subseteq \mathbb K(x)[z]$, and then $S^{-1}\mathcal B=\mathbb K(x)[z]=\mathbb K(x)^{[1]}$. Finally, for item (5), let $a,b \in \mathcal B$ such that $ab \in S=\mathbb K[x]\setminus \{0\}$. Then $ab$ is an unity of $S^{-1}\mathcal B=\mathbb K(x)^{[1]}$ and, thus, $a,b$ are unities of $\mathbb K(x)$. Hence, $a,b \in \mathbb K(x)\cap \mathcal B=\mathbb K[x]$. Therefore $\mathbb K[x]$ is factorially closed in $\mathcal B$ and $\mathcal B^{*}=\mathbb K^{*}$.\endproof

\section{The set \texorpdfstring{$\lnd_{\mathbb K[x]}(\mathcal B)$}{}.}

Let $d$ a derivation on the polynomial ring  $\mathbb K[X,Y,Z]$ such that $$d(X)=0, \quad d(Y)= \varphi_Z (X,Z), \quad \text{ and } \quad d(Z)=f(X).$$ It follows that $d$ is locally nilpotent, since it is easy to see that $d$ is triangular. Moreover, note that $d(F(X,Y,Z))=0$ and, thus, $d$ induces a locally nilpotent derivation $\mathcal D$ on $\mathcal B$ given by 
$$\mathcal D(x)=0, \quad \mathcal D(y)=\varphi_{z}(x,z), \quad \text{ and } \quad \mathcal D(z)=f(x).$$ 

\begin{theorem} \label{lndkx} Let $\mathcal D \in \lnd(\mathcal B)$ as constructed above. Then, 
\begin{itemize}
	\item[(i)]  $\ker \mathcal D= \mathbb K[x]$
	\item[(ii)] $\lnd_{\mathbb K[x]}(\mathcal B)=\{ h \mathcal D \mid h \in \mathbb K[x] \}$;
	\item[(iii)] $\mathcal D$ is irreducible. 
\end{itemize}
\end{theorem}

\proof 
By definition of $\mathcal D$, it is clear that We have $\mathbb K[x] \subseteq \ker \mathcal D \subseteq \mathcal B$. As we have $trdeg_{\mathbb K[x]}(\mathcal B)=1$ and, by Lemma \ref{tjik}(1), $trdeg_{\ker \mathcal D} (\mathcal B)= 1$, it follows that $trdeg_{\mathbb K[x]} (\ker \mathcal D)=0$. So $\ker \mathcal D$ is algebraic over $\mathbb K[x]$. Therefore, from Lemmas \ref{gmn}(5) and \ref{prop}(2), we get $\mathbb K[x] = \ker \mathcal D$.

Let $D \in \lnd_{\mathbb K[x]}(\mathcal B)$, i.e, $D$ is nonzero and $\ker D= \mathbb K[x]$. Since  $\mathcal D(z) =f(x) \in \mathbb K[x]$ we have, by  Lemma \ref{tjik}(3), that $\mathcal D(z)D=D(z)\mathcal D$, where $D(z)\in \mathbb K[x]\backslash\{0\}$. So 
\begin{equation}
f(x)D=D(z)\mathcal D. 
\end{equation}
Now, from Lemma \ref{gmn}(2), we know that $D(y)=\sum_{i<d}a_iz^i$, where $a_i \in \mathbb K[x,y]$, therefore
\begin{equation}
\displaystyle  \sum_{i<d}(f(x)a_i)z^i=f(x)D(y)=D(z)\mathcal D(y)=D(z)\varphi_z(x,z).
\label{isoma}
\end{equation}
As we have $D(z) \in \mathbb K[x]\setminus \{0\}$ and $\varphi_{Z}(X,Z)= dZ^{d-1}+ \sum_{i<d-1}ib_i(X)Z^{i-1}$, it follows that $f(x)a_{d-1}=dD(z)$ and $f(x)a_i=b_iD(z)$, for all $i<d-1$. Therefore, $D(z)=f(x)h$ with $h \in k[x]$. Thus, by Equation (1) we get $D=h\mathcal D$. 

It remains to show that $\mathcal D$ is irreducible. By Lemma \ref{clnd} we have $\mathcal D=hD_0$ for some $h \in k[x]$ and some irreducible $D_0 \in \lnd_{\mathbb K[x]}(\mathcal B)$. We showed above that $D_0=h_0 \mathcal D$ for some $h_0 \in k[x]$, so $\mathcal  D=hh_0 \mathcal D$ and, hence, $h \in \mathbb K^{*}$. Hence, $\mathcal D$ is irreducible. \endproof

\section{The ML-invariant and the Derksen invariant of \texorpdfstring{$\mathcal B$}{\mathcal B}}

We start this section recalling the definition of the $ML$-invariant introduced by Makar-Limanov (some authors refer to this invariant as $AK$-invariant; we follow G. Freudenburg conform \cite{gf1}).

\begin{definition}
Let $A$ be a ring. The $ML$-invariant of $A$, or the ring of absolute constants, is defined as the intersection of the kernels of all locally nilpotent derivation of $A$ and it will be denoted by  $ML(A)$.
\end{definition}

Next is the main result of this Section: 

\begin{theorem}
The $ML$-invariant of $\mathcal B$ is $\mathbb K[x]$. \label{akkx}
\end{theorem}

This Theorem  implies that $\mathbb K[x]\subseteq \ker D$ for all $D \in \lnd(\mathcal B)$ and, by Theorem \ref{lndkx}, there exists 
$D \in \lnd(\mathcal B)$ such that $\ker D=\mathbb K[x]$. In particular,  we have: 

\begin{corollary} \label{lndb} If $\mathcal D$ is the derivation of $\mathcal B$ given by $\mathcal D(x)=0$, $\mathcal D(y)=\varphi_z(x,z)$ and $\mathcal D(z)=f(x)$, then $\lnd(\mathcal B)=\{h \mathcal D \mid h \in\mathbb K[x] \}$. \hfill\qedsymbol
\end{corollary}

Remember that $\mathcal B=\mathbb K[x,y,z]$, where $f(x)y=\varphi(x,z)$, with $r=\deg_Xf(X) >1$ and $d=\deg_Z\varphi(X,Z)>1$. Note that $\mathcal B$ is a subring of $T=\mathbb K[x,f(x)^{-1},z]$, since we have \linebreak $y=f(x)^{-1}\varphi(x,z)\in T$.

Following the strategies of \cite{lml1}, we will construct a so called weight $\mathbb{Z}$-filtration in $T$ and then we use it to induce a $\mathbb Z$-filtration in $\mathcal B$. It will be done because it is possible to obtain some information on a locally nilpotent derivation by passing to its corresponding homogeneous locally nilpotent derivations induced by different filtrations on the ring.

In order to that, first we define a filtration on $R=\mathbb K[x,f(x)^{-1}]$ as follows: for each $n\in \mathbb{Z}$ define the linear subspace $C_n$ of $R$  by setting $$C_n= \begin{cases} \{ax^n \mid a \in \mathbb K\}, & \mbox{ if } n \geq 0, \\
\{a\frac{x^i}{f(x)^{j}}\mid a \in K\}, & \mbox{ if } n=-jr +i, \mbox{ where } 0\leq i \leq r-1 \mbox{ and } j\geq 1.
\end{cases}$$
Notice that, in general, $C_nC_m \nsubseteq  C_{n+m}$ for $n,m\in \mathbb{Z}$. Define $R_k= \bigoplus_{i\leq k}C_i$, for each $k\in \mathbb{Z}$.  We have the following:

\begin{lemma}
\begin{itemize}
	\item[(i)] $R=\displaystyle\bigoplus_{n \in \mathbb{Z}} C_n$.
	\item[(ii)] If $\alpha\in C_n\setminus\{0\}$ and $\beta\in C_m\setminus\{0\}$, then $\alpha\beta\in R_{n+m}\setminus R_{n+m-1}$, for all $n,m\in \mathbb{Z}$. In particular, $R_nR_m \subseteq R_{n+m}$ and $R_n\subset R_{n+1}$ for all  $n,m\in\mathbb{Z}$, $R=\bigcup R_n$, and $\bigcap_n R_n=\{0\}$, that means $\{R_n \mid n \in \mathbb Z\}$ is a proper $\mathbb{Z}$-filtration of $R$.
\end{itemize}
\end{lemma}

\proof (i) Note that every element $\alpha \in R$ is of the form $\alpha=\frac{h(x)}{f(x)^m}$, where $h(x)\in \mathbb K[x]$ e $m\geq 0$. The statement immediate from the following fact: given $h(x)\in \mathbb K[x]\setminus \{0\}$,  by the Euclidean algorithm, there exist $0\leq n_1< \dots <n_s \in \mathbb{Z}$ and $h_1(x), \dots, h_s(x) \in \mathbb K[x]\setminus \{0\}$, with $gr(h_i(x))\leq r-1$, such that $h(x)=\sum_i h_i(x)f(x)^{n_i}$. 

(ii) If $m,n\ge 0$, the claim is immediate. We will treat the case $m\geq 0$ and $n<0$, and the case $m,n<0$ is analogous.
We have $n=-jr +i$, where $j\geq1$ and $0\leq i \leq r-1$. We want to see if $x^m\frac{x^i}{f(x)^{j}} \in R_{n+m}$. If $m+i\leq r-1$ note that $\frac{x^{m+i}}{f(x)^{j}} \in C_{n+m} \subseteq R_{n+m}\setminus R_{n+m-1}$ and the result follows. Now, let us suppose that $m+i\geq r$. We know that $x^{m+i}$ is uniquely written in the form $x^{m+i}=g_1(x)f(x)^{m_1}+ \dots+ g_s(x)f(x)^{m_s}$, where $0\leq m_1<\dots <m_s$ and $\deg(g_{k}(x))\leq r-1$ for all  $1\leq k \leq s$. From this we get $ \frac{x^{m+i}}{f(x)^{j}}=\frac{g_1(x)f(x)^{m_1}}{f(x)^{j}}+ \dots+ \frac{g_s(x)f(x)^{m_s}}{f(x)^{j}}$.  By looking at each term $\frac{g_{k}(x)f(x)^{m_k}}{f(x)^{j}}$, $1\leq k < s$,  in the right hand side of this equality, we have $\deg(g_k)+rm_k \leq  r-1 + rm_k =  (1+m_k)r-1 \leq m_sr-1 <  m_sr  \leq  m_sr + \deg(g_k)$, since $\deg(g_k)<r-1$ and $m_k<m_s$. Thus, $\deg(g_k)+rm_k -jr < \deg(g_s)+rm_s -jr=m+n$, because $\deg(f(x)^{-j}x^{m+i})=\deg(g_s(x)f(x)^{m_s})$. Hence, $\frac{x^{m+i}}{f(x)^{j}} \in R_{m+n}\setminus R_{m+n-1}$.
\endproof

We now introduce a weight $\mathbb{Z}$-filtration in the ring $T=\mathbb K[x, {f(x)}^{-1}, z]$: since each element of $T$ is a finite sum of elements of the form $c_iz^j$, where $j\in \mathbb N$ and $c_i \in C_i$ for $i \in \mathbb{Z}$, to define a such filtration, it suffices to give integer weights $\mu\geq 1$ to $x$ and $\nu$ to $z$. In this case, the weight of a monomial $c_iz^j$ is $i\mu +j\nu$ and the weight of an element $w \in T$ is the maximal weight of the monomials which appear in $w$. By setting $T_n=\text{span}_\mathbb K \{ c_iz^j\mid i\mu +j\nu \leq n\}$ we obtain a $\mathbb Z$-filtration in $T$, which induces one in $\mathcal B$ given by $B_n= T_n\cap \mathcal B$. Additionally, we can extend it to a weight $\mathbb{Z}$-filtration in $\mathbb K(x,z)$ by defining the weight of  $\frac{p}{q} \in \mathbb K(x,z)$ as the difference between the weights of $p$ and $q$. It is well known that the associated graded algebra $Gr(\mathbb K(x,z))$ is isomorphic to the subalgebra of $\mathbb K(x,z)$ consisting of fractions with homogeneous denominators.

\

\noindent
\textbf{The proof of Theorem \ref{akkx}}. Let $D \in \lnd(\mathcal B)$. By Lemma \ref{tjik} parts (1) and (2), there exist $l \in ker(D)\setminus \mathbb K$ and $h \in \mathbb K(x,z)$ such that $D(g)=h J(l,g)$ for all $g\in \mathcal B$, where $J$ is the Jacobian relative to $x$ and $z$. To conclude this proof, it suffices to show that $l \in \mathbb K[x]$ (it will be done in the next lemma). In fact, in this case $D(l)=0$ and, since we have $\ker D$ factorially closed in $\mathcal B$, we conclude that $x\in \ker D$. Therefore, $\mathbb K[x] \subseteq \ker D$ and the proof is completed by Theorem \ref{lndkx}.
 
\begin{lemma} 
$l \in \mathbb K[x]$.
\end{lemma}

\proof First we will show that $l \in \mathbb K[x,z]$. Recall that already have  $l \in ker(D)\setminus \mathbb K$. By Lemma \ref{gmn} we can uniquely write $l=l_m(x,z)y^m + \dots + l_1(x,z)y+l_0(x,z)$, where $m\geq 0$, $l_m(x,z)\neq 0$, and $l_s(x,z)=0$ or $\deg_{z}(l_s(x,z))\leq d-1$  for all $s$ such that $0\leq s \leq m$. 

Suppose that $m\geq 1$. Once we have $\deg_{z}(l_s(x,z))\leq d-1$ and $y=f(x)^{-1}\varphi(x,z)$, we  can assign weight $1$ to $x$ and a sufficiently large weight $\nu$ to $z$ such that $gr(l)=x^kz^tgr(y)^m\in Gr(\mathcal B)$ and $gr(h)=x^cz^d\in Gr(\mathbb K(x,z))$ for some $c,d\in \mathbb{Z}$, after identifying $x$ and $z$ with $gr(x)$ and $gr(z)$, respectively.

Remember that $D(g)=hJ(l,g)$ for $h \in \mathbb K(x,z)$ and that $\mathcal B=\bigcup B_n$, where $\{B_n \}$ is a $\mathbb{Z}$-filtration of $\mathcal B$ induced by the $\mathbb{Z}$-filtration $\{T_n \}$ of $T$. We will see that $D(B_n)\subseteq B_{n+k}$ for all $n \in \mathbb{Z}$ when $k=\deg(l)+ \deg(h)-\deg(x) -\deg(z)$: suppose $x^iz^j\in B_n$, with $i,j\geq 0$, then $D(x^iz^j)=hJ(l,x^iz^j)=h(jx^iz^{j-1}l_x - ix^{i-1}z^jl_z).$
Since $jhx^iz^{j-1}l_x$ and $ihx^{i-1}z^jl_z$ are elements in $B_{n+k}$, we get $D(B_n)\subseteq B_{n+k}$. Now, let suppose $x^if(x)^{-j}z^k\in B_n$, with $i,j,k\geq 0$. Then, $D(x^if(x)^{-j}z^m)= hJ(l,x^if(x)^{-j}z^m)=hJ(l,x^if(x)^{-j}z^m) =h(m\frac{x^i}{f(x)^{j}}z^{m-1}l_x - (i\frac{x^{i-1}}{f(x)^{j}}+j\frac{x^{i}f'(x)}{f(x)^{j+1}})z^ml_z)$, since each summand of this last equality lies in  $B_{n+k}$, we conclude that $D(B_n)\subset B_{n+k}$ for all $n\in \mathbb{Z}$.  

Let $\mathbb{D} \in \lnd(Gr(\mathcal{B}))$ the derivation induced by $D$ for $k=\deg(l)+ \deg(h)-\deg(x) -\deg(z)$ (cf. Section \ref{pre}), i.e., given $g\in B$, 
\[
\mathbb{D}(gr(g))=\begin{cases} gr(D(g)), & \mbox{ if } \deg(D(g))=\deg(g)+k\\ 0,  &\mbox{ if } \deg(D(g))<\deg(g)+k \end{cases}.
\]
Since $D(g)=hJ(l,g)$, we have $\mathbb{D}(gr(g))=gr(h)J(gr(l),gr(g))$. If $gr(l)=x^vz^wgr(y)^m$ we get $v=0=w$, since $\mathbb{D} (gr(l))=0$ and $x,z,gr(y)$ generate the $Gr(\mathcal B)$. Thus, $gr(l)=gr(y)^m$ and $\mathbb{D}(gr(g))=mgr(y)^{m-1}gr(h)J(gr(y),gr(g))$. But, $gr(y)=\frac{z^d}{f(x)}$ then we can write $gr(y)^{m-1}gr(h)=\frac{x^az^b}{f(x)^{(m-1)}}$. 

From Lemma 5 we know that $\mathcal B$ is generated as $\mathbb K$-vector space by $\{x^iz^{k}y^j \mid i,j,k\in \mathbb{N}, 0\leq k\leq d-1\}$ and thus $\{x^iz^{k}gr(y)^j \mid i,j,k\in \mathbb{N}, 0\leq k\leq d-1\}=\{\frac{x^iz^{k+jd}}{f(x)^j} \mid i,j,k\in \mathbb{N}, 0\leq k\leq d-1\}$ is a basis of the $\mathbb K$-vector space $Gr(\mathcal B)$. 

Defining the derivation $\dot{\mathbb{D}} \in \lnd(Gr(\mathcal B))$ by $\dot{\mathbb{D}}(g)=\frac{x^az^b}{f(x)^{(m-1)}}J(gr(y),g)$, where $a,b\in \mathbb{Z}$ and $g\in Gr(\mathcal B)$, we get $\dot{\mathbb{D}}(x)=\frac{-dx^az^{b+d-1}}{f(x)^{m}}$ and $\dot{\mathbb{D}}(z)=\frac{-x^az^{b+d-1}f'(x)}{f(x)^{m+1}}=\frac{-x^az^{b+d}x^{r-1}}{f(x)^{m+1}}$. So $\dot{\mathbb{D}}(x)$ and $\dot{\mathbb{D}}(y)$ will belong to $Gr(\mathcal B)$ if $(a-rm,b+d-1)=(\beta_x-\alpha_xr,\gamma_x+\alpha_xd)$  and $(a+r-1-r(m+1),b+d)=(a-1-rm,b+d-1)=(\beta_z -r\alpha_z,\gamma_z+\alpha_zd)$, where $\alpha_x, \alpha_z, \beta_x, \beta_z, \gamma_x, \gamma_z \in \mathbb{N}$. From these equalities we obtain: $a+(\alpha_x-1)r\geq a+(\alpha_x -m)r \geq 0$, \ $b-1-(\alpha_x-1)d \geq 0$, \ $a-1+(\alpha_z-1)r\geq a-1+(\alpha_z-m)r\geq 0$, \ $b-(\alpha_z-1)d \geq 0$. So $\frac{b-1}{d}\geq \alpha_x-1 \geq \frac{-a}{r}$ and $\frac{b}{d}\geq \alpha_z-1 \geq\frac{1-a}{r}$.

Now using the fact that $\dot{\mathbb{D}}$ is locally nilpotent and denoting ${\deg}_{\dot{\mathbb{D}}}(x)=p$ and ${\deg}_{\dot{\mathbb{D}}}(z)=q$ we have $rp-dq=0$, because $\deg(gr(l))=0$. But $\dot{\mathbb{D}}(x)=-dx^az^bf(x)^{-1}z^{d-1}$ and then $p-1=(a-r)p+(b+d-1)q$. So, $p=d\Delta^{-1}$ and $p=r\Delta^{-1}$, where $\Delta=r(1-b)+d(1-a)$. Hence $\Delta>0$ and divides $r$ and $d$. Note that $\alpha_x-1, \alpha_z-1 \in [\frac{-a}{r}, \frac{b}{d}]$ and the range of this interval is $\frac{b}{d} + \frac{a}{r}=\frac{r+d-\Delta}{rd}<\frac{r+d}{rd}\leq 1$, due to $r>1$ and $d>1$. Hence, $\alpha_x-1=\alpha_z-1$ and $\frac{b-1}{d}\geq \frac{1-a}{r}$, which implies $0\geq \Delta$, that is a contradiction. Therefore, we have $l \in \mathbb K[x,z]$. 

Now, suppose that $l\in \mathbb K[x,z]\setminus \mathbb K[x]$ and $\mu=1$, $\nu=N$ with $N$ sufficiently large such that $gr(l)$ and $gr(h)$ have the same degree in $z$ that $l$ and $h$ have, respectively. 
We saw that $\mathbb{D} (gr(g))=gr(h)J(gr(l),gr(g))$. 
Since  $l \not\in \mathbb K[x]$, we have $gr(l)=\alpha x^iz^j$ with $j>0$. So, $z$ divides $gr(l)$ and, hence, $z\in \ker(\dot{\mathbb{D}})$. 
Finally, as $f(x)gr(y)= z^d$, we get $gr(y) \in ker(\dot{\mathbb{D}})$ and we conclude that $\dot{\mathbb{D}}$ is identically zero. However, $\dot{\mathbb{D}}(x)=jx^iz^{j-1}\neq 0$, which is a contradiction. Hence, $l \in \mathbb K[x]$, as we desired to proof. \endproof

To conclude this section, we recall the definition of the invariant introduced by Derksen (see \cite{derk}).

\begin{definition}
Let $A$ be a $\mathbb K$-algebra. The Derksen invariant of $A$ is the $\mathbb K$-algebra generated by the union of the kernels of all nontrivial locally nilpotent derivation of $A$ and it will be denoted by  $HD(A)$.
\end{definition}

\begin{theorem} \label{HD}
The $HD$-invariant of $\mathcal B$ is $\mathbb K[x]$. 
\end{theorem}

\proof It is a consequence of Corollary \ref{lndb} where we get $\ker D=\mathbb K[x]$ for all $D\in \lnd(\mathcal B)$.
\endproof

\section{On automorphisms of \texorpdfstring{$\mathbb K[X,Y,Z]/(f(X)Y - \varphi(Z))$}{}} 

We shall describe a set of generators for the group of $\mathbb K$-automorphisms of the ring $\mathcal B=\mathbb K[X,Y,Z]/(f(X)Y - \varphi(Z))$, where $\varphi(Z)\in \mathbb K[Z]$ has degree $d>1$, denoted by $\aut(\mathcal B)$. There is a complete description of $\aut(\mathcal B)$ when $f(X)=X^r$, with $r\ge 1$, in \cite{dda1,lml3}. For this reason, we assume that $f(X)$ has at least one nonzero root. 

As we are interested in $\mathbb K$-automorphisms of $\mathcal B$, notice that we can assume, without loss of generality, that the coefficient of the term with degree $d-1$ of $\varphi(Z)$ and the coefficient of the term with degree $r-1$ of $f(X)$ are equals to zero (to see that, just change $Z$ by $Z-a$ and $X$ by $X-b$ for a convenient choice of $a,b\in \mathbb K$).

In order to obtain this characterization, it will be clear that is essential to know that the $ML$-invariant of $\mathcal B$ is $\mathbb K[x]$.

\begin{lemma} Let $g(X)$ be a monic polynomial of degree $r\geq 2$ in $\mathbb K[X]$ with at least one nonzero root and such that the coefficient of the term of degree $r-1$ is zero. If $\lambda,\beta \in \mathbb K$ and $\lambda\neq0$, then $\lambda^rg(X)=g(\lambda X+\beta)$ if, and only if, $\beta=0$, $\lambda$ is a $s$-th root of unity and $g(X)=X^ih(X^s)$ with $h(X) \in \mathbb K[X]$ monic. \label{xaxzaz}
\end{lemma}

\proof Suppose that $\lambda^rg(X)=g(\lambda X+\beta)$. By taking the Taylor expansion at $X=0$ for $g(X)$ we have
\begin{equation}
g(X)= \frac{g^{(r)}(0)}{r!}X^r + \frac{g^{(r-1)}(0)}{(r-1)!}X^{r-1}+\dots+ g'(0)X +g(0).
\label{taylor}
\end{equation}
Now, deriving successively the expression $\lambda^rg(X)=g(\lambda X+\beta)$, we get
\begin{equation}
\lambda^rg^{(n)}(X)= \lambda^ng^{(n)}(\lambda X+ \beta). \label{idfn}
\end{equation}
and, by calculating the $(r-1)$-th derivative of Equation \eqref{taylor}, we have 
\begin{equation}
g^{(r-1)}(X)= g^{(r)}(0)X + g^{(r-1)}(0). \label{fr1}
\end{equation}
From \ref{idfn} and \ref{fr1}, $\lambda g^{(r-1)}(0)= g^{(r-1)}(\beta)= g^{(r)}(0)\beta + g^{(r-1)}(0)$. But, by hypothesis, $g^{(r-1)}(0)=0$ and $g^{(r)}(0)=1$, thus $\beta=0$ and then $\lambda^rg(X)=g(\lambda X)$.

Let us now show that $g(X)=X^ih(X^s)$ as stated. By hypothesis, $g(X)$ has a nonzero root, so there exist $t\in \mathbb{N}$  and positive integers $0\leq n_0<\ldots<n_t\leq r-2$ such that $g(X)=X^r+\sum_{l=0}^{t} a_lX^{n_l}$, where $a_0,a_1,\ldots, a_{t} \in \mathbb K^*$. Since $\lambda^rg(X)=g(\lambda X)$ we have $\lambda^ra_{i}= \lambda^{n_i}a_{i}$ and thus $\lambda^{r-n_i}=1$ for $i=0,1,\ldots,t$. Therefore $\lambda$ is a $(r-n_i)$-th root of unity. Let $s$ be the order of  $\lambda$, so $s$ divides $mdc\{r-n_0, \ldots, r-n_t\}$. Writting  $r=ns+i$, for some $n \in \mathbb Z_+$ and $0\leq i<s$, we have $\lambda^ig(X)=\lambda^rg(X)=g(\lambda X)$ and then $\lambda^ia_l=a_l\lambda^{n_l}$ for $l=0,1,\ldots,t$, which implies $\lambda^{n_l-i}=1$ for $l=0,1, \ldots,t$. Therefore, there exist $m_l\in\mathbb Z$ such that $n_l=sm_l+i$ for all $l\in \{1, \ldots, t\}$. So 
$g(X)=X^r+\sum_{l=0}^{t} a_lX^{n_l}=X^{sn+i} + \sum_{l=0}^{t} a_iX^{sm_l+i} = X^i(X^{sn} + \sum_{l=0}^{r-2} a_lX^{sm_l})$.
Hence, $g(X)=X^ih(X^s)$ with $h(X)=X^{sn} + \sum_{l=0}^{t} a_iX^{m_l}$.\endproof

We now start the study of $\mathbb K$-automorphisms of $\mathcal B$. Let $T$ be an $\mathbb K$-automorphism of $\mathcal B$. First, note that $T(\mathbb K[x])=\mathbb K[x]$, since $\mathbb K[x]$ is the $ML$-invariant of $\mathcal B$ (see \cite{gf1}). Thus, $T(x)=\lambda x + \beta$ for some $\lambda,\beta \in \mathbb K$ and $\lambda \neq0$. Let
$\mathcal D=f(x)\frac{\partial}{\partial z}$ be the locally nilpotent derivation of $\mathcal B$ conform Corollary \ref{lndb}. Then  $\mathcal{D}^2(z)=0$, and, again by Corollary \ref{lndb}, for any locally nilpotent derivation $\partial$ of $\mathcal B$ we have $\partial^2(z)=0$. Therefore, ${T}^{-1} \partial^2(T(z))=0$ for any $\partial \in \lnd(\mathcal B)$ and then $\partial^2(T(z))=0$ for any $\partial \in \lnd(B)$.  
Now, as we have $\mathcal D(T(z))=g(x,y,z)$, by taking $\partial=f(x)^m \mathcal D$ with $m$ large enough such that $\partial(T(z))= h(x,z)$. It follows that  $T(z)=\alpha (x)z +b(x)$, with $\alpha (x)\neq 0$, since $\partial (z)\in \mathbb K[x]$, $\mathbb K[x]\subseteq \ker(\partial)$, and $x,z$  are algebraically independent. Since $T$ has an inverse, we conclude that $\alpha (x) \in \mathbb K^*$. Hence, $T(z)=cz +b(x)$ with $c \in \mathbb K^*$.

\begin{lemma} Let $T$ be an $\mathbb K$-automorphism of $\mathcal B$. Then $T(z)=\alpha z+b(x)$ and $T(x)=\lambda x$, with $\lambda, \alpha \in \mathbb K^*$, $b(x) \in \mathbb K[x]$, $b(x)\equiv 0 \mod f(x)$, and $\varphi(\alpha z)=\alpha^d\varphi( z)$. Moreover, $f(x)=x^ih(x^s)$, with $h\in\mathbb K^{[1]}$, and $\lambda^s=1$. \label{lautb}
\end{lemma}

\proof Let $T$ be an $\mathbb K$-automorphism of  $\mathcal B$ and $\mathcal D=f(x)\partial/\partial z$ the locally nilpotent derivation given by Corollary \ref{lndb}. We know, 	by the previous remark, that $T(x)=\lambda x + \beta$ with $\lambda,\beta \in \mathbb K$ e $ \lambda \neq0$, $T(z)=cz +b(x)$ where $c \in \mathbb K$ and $b(x)\in \mathbb K[x]$. Note that
\begin{equation}
T \mathcal DT^{-1}(z)=T \mathcal D(c^{-1}z-T^{-1}(b(x))=c^{-1}T(f(x))=c^{-1}f(\lambda x + \beta). \label{tdt}
\end{equation}
By Corollary \ref{lndb}, $f(x)$ divides $D(z)$ for any $D\in LND(\mathcal B)$. Then, by Equation \eqref{tdt}, $f(x)$ divides $f(\lambda x +\beta)$, because $T\mathcal D{T}^{-1} \in LND(\mathcal B)$, i.e., $\delta f(x)=f(\lambda x +\beta)$ with $\delta \in \mathbb K^*$. Therefore, we find that $\delta=\lambda^r$. Thus, $\lambda^rf(x)=f(\lambda x +\beta)$ and, by Lemma \ref{xaxzaz},	we get the result  for $T(x)$ and the polynomial $f(x)$.
	
Applying $T$ to the equality $f(x)y=\varphi(z)$, we have $\lambda^rf(x)T(y)=f(\lambda x)T(y)=T(f(x))T(y)=T(\varphi(z))=\varphi(cz +b(x))$.
Note that $\varphi(cz +b(x))=c^d\varphi(z)+\nu(x,z)$, where $\nu(x,z) \in k[x,z]$ and $\deg_z(\nu)\leq d-1$. Therefore we have $\lambda^{r}T(y)=f(x)^{-1}(c^d\varphi(z)+\nu(x,z))=c^dy+\nu(x,z)f(x)^{-1}$, that means $\nu(x,z)f(x)^{-1} \in \mathbb K[x,y,z]$ and  then $\nu(x,z) \equiv0 \mod f(x)$. Now, note that $0\equiv_{f(x)} \nu(x,z)= \varphi(cz +b(x))-c^d\varphi(z)=c^{d-1}z^{d-1}b(x) + \sigma$, where $\sigma \in k[x,z]$ and $\deg_z(\sigma)<d-1$. We conclude, by Lemma \ref{gmn}, that $b(x)\equiv0 \mod f(x)$. As $\nu(x,z)=\varphi(cz +b(x))-c^d\varphi(z)$ we have $\varphi(cz +b(x))\equiv\varphi(cz) \mod f(x)$ and since $b(x)\equiv0 \mod f(x)$ we obtain that $\varphi(cz)\equiv c^d\varphi(z) \mod f(x)$, but this is only possible when $\varphi(cz)-c^d\varphi(z)=0$.
\endproof

\begin{theorem} \label{autb} The group $\aut(\mathcal B)$ is generated by the following $\mathbb K$-automorphisms: 
\begin{enumerate}
\item $H$ defined by $H(x)=x$, $H(y)= y + [\varphi(z+h(x)f(x)) -\varphi(z)]f(x)^{-1}, \text{ and } H(z)=z + h(x)f(x), \text{ where } h \in \mathbb K^{[1]}$,
\item $T$ defined by $T(x)=\lambda x$, where  $\lambda^s=1$, $T(y)=\lambda^jy$ and $T(z)=z$, if $f(X)=X^jh(X^s)$, with  $h \in \mathbb K^{[1]}$ having a nonzero root,
\item $R$ defined by $R(x)=x$, $R(y)=\lambda^dy$ and $R(z)=\lambda z$, if $\varphi(z)=z^d$,
\end{enumerate}
and
\begin{enumerate}
\setcounter{enumi}{3}
\item $S$ defined by  $S(x)=x$, $S(y)=\mu^iy$ and $S(z)=\mu z$ where $\mu^m=1$, if $\varphi(z)=z^i\varphi(z^m)$, with $\phi\in \mathbb K^{[1]}$.
\end{enumerate}

\end{theorem}

\proof It is clear that all of these maps are $\mathbb K$-automorphisms of $\mathcal B$. It follows from Lemma \ref{lautb} that any $\mathbb K$-automorphism is a composition of an $\mathbb K$-automorphism  $W$, given by $W(x)=\lambda x$ and $W(z)=z$, an $\mathbb K$-automorphism $Y$, given by $Y(x)=x$ and $Y(z)=\alpha z$, and the $\mathbb K$-automorphism $H$ as above. The $\mathbb K$-automorphisms of $(3)$ and $(4)$ treat the case $\varphi(z)$ when $\varphi(\alpha z)=\alpha^d\varphi(z)$, with $\alpha\neq1$. The $\mathbb K$-automorphisms in all other cases deal with all possible cases involving the $\mathbb K$-automorphism $W$ with $\lambda\neq1$. \endproof

\begin{remark}
For the particular case where $f(x)= x^ih(x)$, with  $h(x)=x^{r-i}+\sum_{l=0}^{t} a_lx^{m_l},$ where $a_l\in \mathbb K^{*}$, for $0\leq l\leq t$, $0=m_0<m_1<\dots <m_t<r-i$ and $\text{mdc}\{r-i,m_1,\dots, m_t\}=1$, we have that every $\mathbb K$-automorphism of $\mathcal B$ is a $\mathbb K[x]$-automorphism. Additionally, if we have $\varphi(z)=z^{d}+b_{d-2}(x)z^{d-2}+\dots+b_0(x)\in \mathbb K[x,z]$, then the group of $\mathbb K$-automorphisms of $\mathcal B$ is generated by the group of $\mathbb K$-automorphisms of items $(1)$, $(3)$, and $(4)$ of the previous theorem.
\end{remark}

\begin{example} Let $ \mathcal B=\mathbb{C}[x,y,z]$ such that $f(x)y=\varphi(z)$, where
 $\varphi(Z)=Z^3+Z+1 \in \mathbb{C}[Z] $, $f(X)=X^2h(X^4)$, and $h(X)=X^5+2X^4+X^2-2 \in \mathbb{C}[X]$. By Theorem  \ref{autb}, we have an $\mathbb K[x]$-automorphism $H$ of $\mathcal B$ defined  by
\[
H(x)=x,\qquad H(z)=z +(x^2+1)f(x)= x^{24}+x^{22}+2x^{20}+2x^{18}+x^{12}+x^{10}-2x^{4}-2x^{2},
\]
\[
\begin{array}{rcl} 
H(y)&=&y+[\varphi(z+(x^2+1)f(x))-\varphi(z)]{f(x)}^{-1} \\
    &=& y+(1 + x^2)(1 + 4 x^4 + 8x^6 + 4x^8 - 4x^{12} - 8x^{14} - 4x^{16}\\
		& &-7x^{20} - 14x^{22} - 11 x^{24} - 8x^{26} + 8 x^{30} + 6x^{32} + 4x^{34} \\
		& &+6 x^{36} + 8 x^{38} + 8 x^{40} + 8 x^{42} + 5 x^{44} + 2 x^{46} + x^{48}-6x^2 z - 6x^4z  \\
		& &+3x^{10}z + 3x^{12}z + 6x^{18}z + 6x^{20}z +3x^{22}z + 3x^{24}z + 3z^2). 
\end{array}
\]

\end{example}


\begin{thebibliography}{99}

\bibitem{cra1} 
Anthony J. Crachiola, \textit{ On automorphisms of Danielewski surfaces}, J. Algebraic Geom. {\bf 15} (2006), no. 1, 111-132.

\bibitem{dda1} 
D. Daigle, \textit{On locally nilpotent derivations of $\mathbb K[X_1,X_2,Y]/(\varphi(Y)- X_1X_2)$}, J. Pure and Appl. Algebra
\textbf{181} (2003), 181--208 .

\bibitem{dda2}
D. Daigle, \textit{Locally nilpotent derivations and Danielewski surfaces}, Osaka J. Math.  {\bf 41} (2004),  no. 1, 37--80. 

\bibitem{dda3}
D. Daigle, \textit{Locally nilpotent derivations}, Lecture notes for the September School of algebraic geometry, Luk\c ecin, Poland,
Setember 2003, Avaible at \verb"http://aix1.uottawa.ca/~ddaigle". 

\bibitem{dan}
W. Danielewski, \textit{On the cancellation problem and automorphism groups of affine algebraic varieties}, preprint, (1989)


\bibitem{derk}
H. Derksen, Constructive Invariant Theory and the Linearization Problem, PhD. thesis, University of Basel, 1997

\bibitem{dubo}
A. Dobouloz, \textit{ On a class Danielewski surfaces in affine 3-space}, J. Algebra \textbf{321} (2009), 1797--1812


\bibitem{gf1}
G. Freudenburg, Algebraic Theory of Locally Nilpotent Derivations, Encyclopaedia of Mathematical Sciences, Springer-Verlag (2006).


\bibitem{lml1}
L. Makar-Limanov, \textit{On groups of automorphisms of a class of surfaces}, Israel J. of Math. \textbf{68} (1990), 250--256.

\bibitem{lml2}
L. Makar-Limanov, \textit{On the hypersurface $x + x^2y +z^2 +t^3=0$ in $C^4$ or a $C^3$-like three fold which is not $C^3$}, Israel J. of
Math. \textbf{96} (1996), 419--429 .


\bibitem{lml3}
L. Makar-Limanov, \textit{On the groups of automorphisms of a surface $x^ny=P(z)$}, Israel J. of Math. \textbf{121} (2001), 113--123.


\end{thebibliography}
\end{document}